\documentclass{amsart}
\usepackage{amssymb,amsfonts,latexsym}
\newtheorem{theorem}{Theorem}[section]
\newtheorem{lemma}[theorem]{Lemma}

\newtheorem{corollary}[theorem]{Corollary}
\theoremstyle{definition}

\newtheorem{remark}[theorem]{Remark}

\newcommand{\id}{\text{id}}

\newcommand{\Ad}{\text{Ad}}

\newcommand{\gr}{\text{gr}}

\newcommand{\g}{\mathfrak{g}}

\newcommand{\ot}{\otimes}

\newcommand{\ben}{\begin{enumerate}}
\newcommand{\een}{\end{enumerate}}

\newcommand{\Z}{{\mathbb Z}}

\newcommand{\C}{{\mathbb C}}

\begin{document}
\title[on radically graded finite dimensional Quasi-Hopf algebras]
{on radically graded finite dimensional Quasi-Hopf algebras}
\author{Pavel Etingof}
\address{Department of Mathematics, Massachusetts Institute of Technology,
Cambridge, MA 02139, USA} \email{etingof@math.mit.edu}
\author{Shlomo Gelaki}
\address{Department of Mathematics, Technion-Israel Institute of
Technology, Haifa 32000, Israel}
\email{gelaki@math.technion.ac.il}
\date{January 13, 2005}
\begin{abstract}
In this paper we continue the structure theory of finite dimensional
quasi-Hopf algebras started in \cite{eg} and \cite{g}. First, we
completely describe the class of radically graded finite dimensional
quasi-Hopf algebras over $\mathbb C$, whose radical has prime
codimension. As a corollary we obtain that if $p>2$ is a prime then
any finite tensor category over $\C$ with exactly $p$ simple objects
which are all invertible must have Frobenius-Perron dimension $p^N$,
$N=1,2,3,4,5$ or $7$. Second, we construct new examples of finite
dimensional quasi-Hopf algebras which are not twist equivalent to a
Hopf algebra. For instance, to every finite dimensional simple Lie
algebra ${\mathfrak g}$ and an odd integer $n$, coprime to $3$ if
$\g=G_2$, we attach a quasi-Hopf algebra of dimension $n^{{\rm
dim}{\mathfrak g}}$.
\end{abstract}
\maketitle
\section{Introduction}
In \cite{eo} it is proved that any finite tensor category over
$\mathbb{C}$ with integer Frobenius-Perron dimensions of objects is
equivalent to a representation category of a finite dimensional
quasi-Hopf algebra (the Frobenius-Perron dimension of a
representation coincides with its dimension as a vector space).
Therefore the classification of finite tensor categories with
integer Frobenius-Perron dimensions of objects is equivalent to the
classification of complex finite dimensional quasi-Hopf algebras.
The simplest finite tensor categories to try to understand are those
which have only 1-dimensional simple objects which form a cyclic
group of prime order under tensor product. Equivalently, one is led
to the problem of classifying finite dimensional quasi-Hopf algebras
with (Jacobson) radical of prime codimension.

Let $p$ be a prime, and let $RG(p)$ denote the class of radically
graded finite dimensional quasi-Hopf algebras over $\mathbb C$,
whose radical has codimension $p$. It was shown in \cite{eg} that
any $H\in RG(2)$ is equivalent to a Nichols Hopf algebra $H_{2^n}$,
$n\ge 1$ \cite{n}, or to a lifting of one of the four special
quasi-Hopf algebras $H(2)$, $H_+(8)$, $H_-(8)$, $H(32)$ of
dimensions 2, 8, 8, and 32 (the algebra $H(2)$ is the group algebra
of $\mathbb Z_2$ with a nontrivial associator).

Later, it was shown in \cite{g} that if $H\in RG(p)$, $p>2$, has a
nontrivial associator and if the rank of $H[1]$ over $H[0]$ is $1$,
then $H$ is equivalent to one of the quasi-Hopf algebras $A(q)$ of
dimension $p^3$, introduced in \cite{g}. More precisely, the result
of \cite{g} is formulated under the assumption that $H$ is basic
(i.e., $H/{\rm Rad}(H)=\C[\Z_p]$ with some associator), but by
\cite{eno}, Corollary 8.31, this is automatic.

The purpose of this paper is to continue the structure theory of
finite dimensional quasi-Hopf algebras started in \cite{eg} and
\cite{g}. More specifically, we completely describe the class
$RG(p)$, and construct new examples of finite dimensional quasi-Hopf
algebras which are not twist equivalent to a Hopf algebra.

The structure of the paper is as follows. In Section 2 we recall the
definition of the quasi-Hopf algebras $A(q)$ and $H_{\pm}(p)$.

In Section 3 we show that if $H\in RG(p)$ has a nontrivial
associator, then the rank of $H[1]$ over $H[0]$ is $\le 1$. This
yields the following classification of $H\in RG(p)$, $p>2$, up to
twist equivalence.

(a) Duals of pointed Hopf algebras with $p$ grouplike elements,
classified in \cite{as}, Theorem 1.3.

(b) Group algebra of $\Z_p$ with associator defined by a 3-cocycle.

(c) The algebras $A(q)$.

This result implies, in particular, that if $p>2$ is a prime then
any finite tensor category over $\C$ with exactly $p$ simple objects
which are all invertible must have Frobenius-Perron dimension $p^N$,
$N=1,2,3,4,5$ or $7$.

In Section 4 we construct new examples of finite dimensional
quasi-Hopf algebras $H$, which are not twist equivalent to a Hopf
algebra. They are radically graded, and $H/{\rm Rad}(H)=\C[\Z_n^m]$,
with a nontrivial associator. For instance, to every finite
dimensional simple Lie algebra ${\mathfrak g}$ and an odd integer
$n$, coprime to $3$ if $\g=G_2$, we attach a quasi-Hopf algebra of
dimension $n^{{\rm dim}{\mathfrak g}}$.

{\bf Acknowledgments.} The research of the first author was
partially supported by the NSF grant DMS-9988796.
The second author was supported by Technion V.P.R.
Fund - Dent Charitable Trust - Non Military Research Fund, and by The Israel
Science Foundation (grant No. 70/02-1). He also thanks MIT for its
warm hospitality.
Both authors were supported by BSF grant No. 2002040.
\section{Preliminaries}
All constructions in this paper are done over the field of complex
numbers $\mathbb{C}$.

We refer the reader to \cite{d} for the definition of a quasi-Hopf
algebra and a twist of a quasi-Hopf algebra.

\subsection{} We recall the theory of the radical filtration for
finite dimensional quasi-Hopf
algebras, discussed in \cite{eg}.
It is completely parallel to the classical theory of such
filtration in finite dimensional Hopf algebras.

Let $H$ be a finite dimensional quasi-Hopf
algebra, and $I$ be the radical of $H$.
Assume that $I$ is a quasi-Hopf ideal, i.e.,
$\Delta(I)\subseteq H\ot I + I\ot H$. In categorical terms, this
means that the category of representations ${\rm Rep}(H)$ has
Chevalley property, i.e., the tensor product of irreducible
$H$-modules is completely reducible. This is satisfied, for example,
if $H$ is basic, i.e., every irreducible $H$-module is
1-dimensional.

In this situation, the filtration of $H$ by powers of $I$
is a quasi-Hopf algebra filtration. Thus the associated graded
algebra ${\rm gr}(H)$ of $H$ under this filtration has a natural
structure of a quasi-Hopf algebra.

Let now $\overline{H}$ be a finite dimensional quasi-Hopf algebra
with a $\Z_+$-grading, i.e., $\overline{H}=\oplus_{m\ge 0}\overline{H}[m]$,
with all structure maps being of degree zero.
In this case, $\overline{H}[0]$ is a quasi-Hopf algebra,
$\overline{H}[i]$ is a free module over $\overline{H}[0]$
for all $i$ (by Schauenburg's theorem \cite{s}), and
the radical $\overline{I}$ of $\overline{H}$
is a quasi-Hopf ideal.

One says that $\overline{H}$ is radically graded if
$\overline{I}^k=\oplus_{m\ge k}\overline{H}[m]$, for $k\ge 1$. In this case,
$\overline{H}[0]$ is semisimple, and $\overline{H}$
is generated by $\overline{H}[0]$ and $\overline{H}[1]$.

An example of a radically graded
quasi-Hopf algebra is the algebra ${\rm gr}(H)$ defined above.
Moreover, $H$ is radically graded if and only if ${\rm gr}(H)=H$.

Finally, we observe that if $H$ is radically graded and basic,
then $H[0]={\rm Fun}(G)$ for a finite group $G$, and
the associator (being of degree zero) corresponds
to a class in $H^3(G,\C^*)$.

\subsection{} The following are the simplest
examples of quasi-Hopf algebras not twist equivalent to a Hopf
algebra.

Let $p>2$ be a prime, and $\varepsilon=e^{2\pi {\rm i}/p}$.
If $z\in \Z$, we denote by $z'$ the projection of $z$
to $\Z_{p}$.

Let $s$ be an integer such that $1\le s\le p-1$.
Let $Q=\varepsilon^{-s}$.
The $p-$dimensional quasi-Hopf algebra $H(p,s)$,
is generated by a grouplike element $a$ such that $a^p=1$, with
non-trivial associator
\begin{equation}\label{assoc}
\Phi_s:=\sum_{i,j,k=0}^{p-1}Q^{\frac{-i(j+k-(j+k)')}{p}}{\bf 1}_i\ot
{\bf 1}_j\ot {\bf 1}_k
\end{equation}
where $\{{\bf 1}_i|0\le i\le p-1\}$ is the
set of primitive idempotents of $\Z_p$
(i.e, ${\bf 1}_ia=Q^i{\bf 1}_i$), distinguished elements
$\alpha=a$, $\beta=1$, and antipode $S(a)=a^{-1}$.

Let $s_0\in \Z_p$ be a quadratic nonresidue. It can be shown (by
considering automorphisms $a\mapsto a^m$) that for any $s$,
$H(p,s)$ is isomorphic to $H_+(p):=H(p,1)$ if $s$ is a quadratic
residue, and to $H_-(p):=H(p,s_0)$ if $s$ is a non-quadratic
residue. On the other hand, $H_+(p)$ and $H_-(p)$ are not
equivalent.

Thus it follows from
\cite{eno}, Corollary 8.31, that any $p-$dimensional semisimple
quasi-Hopf algebra is twist equivalent either to $\C[\Z_p]$
or to $H_\pm(p)$.

\subsection{}
The following are examples of $p^3-$dimensional basic quasi-Hopf
algebras with radical of codimension $p$, which are not twist equivalent
to a Hopf algebra.

\begin{theorem} \cite{g}\label{propp3}
Let $p$ be a prime number.

(i) There exist $p^3-$dimensional quasi-Hopf algebras $A(q)$,
parametrized by primitive roots of unity $q$ of order $p^2$, which
have the following structure. As algebras $A(q)$ are generated by
$a,x$ with the relations $ax=q^p xa$, $a^p=1$, $x^{p^2}=0$. The
element $a$ is grouplike, while the coproduct of $x$ is given by the
formula
\begin{equation*}
\Delta(x)=x\ot \sum_{y=0}^{p-1}q ^y {\bf 1}_y + 1\ot (1-{\bf
1}_{0})x + a^{-1}\ot {\bf 1}_{0}x,
\end{equation*}
where $\{\mathbf{1}_i|0\le
i\le p-1\}$ is the set of primitive idempotents of $\C[a]$
defined by the condition $a{\bf 1}_i=q ^{pi}{\bf 1}_i$, the associator
is $\Phi_s$ (where $s$ is defined by
the equation $\varepsilon^{-s}=q^p$), the distinguished
elements are $\alpha=a$,
$\beta=1$, and the antipode is $S(a)=a^{-1}$,
$S(x)=-x\sum_{z=0}^{p-1}q ^{p-z}{\bf 1}_z$.

(ii) The quasi-Hopf algebras $A(q)$
are pairwise non-equivalent. Any finite dimensional radically graded
basic quasi-Hopf algebra $H$ with radical of codimension $p$ and
nontrivial associator, such that
$H[1]$ is a free module of rank $1$ over $H[0]$, is equivalent
to $A(q)$ for some $q$.
\end{theorem}

\section{quasi-Hopf algebras with radical of prime
codimension}
\subsection{The main result} Let $p>2$ be a prime number.
Our main result in this section is the following theorem.

\begin{theorem}\label{rank} Let
$H$ be a radically graded basic quasi-Hopf algebra
with radical of codimension $p$. If the associator of
$H$ is nontrivial, then the rank of $H[1]$
over $H[0]$ is $\le 1$.
\end{theorem}

Theorem \ref{rank} is proved in the next subsection.

Theorem \ref{rank} and the results cited above
imply the following classification result.

\begin{theorem}\label{class} Let $H$ be a radically graded
finite dimensional quasi-Hopf algebra with radical
of codimension $p$. Then $H$ is one of the following quasi-Hopf
algebras, up to twist equivalence:

(a) Duals of pointed Hopf algebras with $p$ grouplike elements,
classified in \cite{as}, Theorem 1.3 (including the group algebra
$\C[\Z_p]$).

(b) The algebras $H_+(p)$ and $H_-(p)$.

(c) The algebras $A(q)$.
\end{theorem}

\begin{proof}
By Corollary 8.31
of \cite{eno}, $H$ is necessarily basic.

If the associator of $H$ is trivial, then
we may assume that $H$ is a Hopf algebra. Thus
$H^*$ is a coradically graded pointed Hopf algebra
with $G(H^*)=\Z_p$. Such algebras are classified in
\cite{as}, Theorem 1.3, so we are in case (a).

If the associator is nontrivial, then
by Theorem \ref{rank}, the rank of $H[1]$ over $H[0]$
is at most 1. If the rank is 0, we are in case
(b). If the rank is 1, we are in case (c)
by Theorem \ref{propp3}.
\end{proof}

We refer the reader to \cite{eo}, for the definition of
a finite tensor category and the notion of its Frobenius-
Perron dimension.

\begin{corollary} Let $p>2$ be a prime.
Let $\mathcal C$ be a finite tensor category,
which has exactly $p$ simple objects which are all invertible.
Then the possible values of the Frobenius-Perron
dimension of $\mathcal C$ are $p^N$, $N=1,2,3$ (for all $p$),
$4$ (for $p=3$ and $p=3k+1$), $5$ (for $p=3$ and $p=4k+1$) and
$7$ (for $p=3$ and $p=3k+1$).
\end{corollary}

\begin{proof} It is clear that the Frobenius-Perron dimension
of objects in $\mathcal C$ are integers. Hence by \cite{eo}, there
exists a quasi-Hopf algebra $A$ such that $\mathcal C={\rm
Rep}(A)$. This quasi-Hopf algebra is basic, so its radical is a
quasi-Hopf ideal and hence $A$ admits a radical filtration. Let
$H:={\rm gr}(A)$ (with respect to this filtration). Then Theorem
\ref{class} applies to $H$, hence the result.
\end{proof}

\subsection{Proof of Theorem \ref{rank}}
Let us assume that $H[1]$ has rank $>1$ over $H[0]$. From this we
will derive a contradiction. We may assume that $H$ has the minimal
possible dimension.

Let $a$ be a generator of $\Z_p$. We have $H[0]=\C[\Z_p]$ with
associator $\Phi_s$ for some $s$.

Let us decompose $H[1]$ into a direct sum of eigenspaces of $a$:
$H[1]=\bigoplus_{r=0}^{p-1} H_r[1]$, where $H_r[1]$ is the space of
$x\in H[1]$ such that $axa^{-1}=Q^rx$ (we recall that
$Q:=\varepsilon^{-s}$). Note that $\mathbf{1}_ix=x\mathbf{1}_{i-r}$
for $x\in H_r[1]$. Also, by Theorem 2.17 in \cite{eo}, $H_0[1]= 0$.

Let $\tilde{H}$ be the free algebra generated by $H[1]$ as a
bimodule over $H[0]$; i.e., $\tilde{H}$ is the tensor algebra of
$H[1]$ over $H[0]$. Then $\tilde{H}$ is (an infinite dimensional)
quasi-Hopf algebra, and we have a surjective homomorphism
$\varphi:\tilde{H}\to H$ (it is surjective since $H$ is radically
graded and hence generated by $H[0]$ and $H[1]$).

Let $q$ be a number such that $q^p=Q$. Define an automorphism
$\gamma$ of $\tilde{H}$ by the formula $\gamma_{|H[0]}=1$ and
$\gamma_{|H_r[1]}=q^r$. (It is well defined since $\tilde{H}$ is
free.)

Let $L$ be the sum of all quasi-Hopf ideals in $\tilde{H}$ contained
in $\bigoplus_{d\ge 2}\tilde{H}[d]$. Clearly, $Ker\varphi\subseteq
L$, so $H$ projects onto $\tilde{H}/L$. However, since $H$ has the
smallest dimension, it follows that $\tilde{H}/L= H$.

Now, $\gamma(L)=L$, so $\gamma$ acts on $H$. Let us define a new
quasi-Hopf algebra $\hat{H}$ generated by $H$ and a grouplike
element $g$ with relations $g^p=a$, $gzg^{-1}=\gamma(z)$ for $z\in
H$. Clearly, $\Ad(a)=\gamma^p$, and $g$ generates a group isomorphic
to $\Z_{p^2}$.

Let $J:=\sum_{i,j}c(i,j)1_i\ot 1_j$, $c(i,j):=q^{-i(j-j')}$, where
$j'$ denotes the remainder of division of $j$ by $p$, be the twist
in $\C[\Z_{p^2}]^{\ot 2}$ defined in \cite{g}. Define $\bar H$ to be
the twist of $\hat{H}$ by $J^{-1}$: $\bar{H}:=\hat{H}^{J^{-1}}$.
Since by \cite{g}, $dJ=\Phi_s$, $\bar H$ is a finite dimensional
radically graded {\bf Hopf} algebra. Since the rank of $\bar H[1]$
over $\bar H[0]$ is $>1$, we have at least $2$ independent over
$\bar H[0]$ skew primitive elements $x_1,x_2\in \bar{H}[1]$ which
are eigenvectors for $\Ad(g)$:
$$gx_1g^{-1}=q^{d_1}x_1,\;\Delta(x_1)=x_1\ot g^{b_1}+1\ot x_1$$
and
$$gx_2g^{-1}=q^{d_2}x_1,\;\Delta(x_2)=x_2\ot g^{b_2}+1\ot x_2.$$
Since $H_0[1]=0$, $d_1,d_2$ must be relatively prime to $p$.
Also, since $H$ has minimal dimension, the algebra $\bar H$ is
generated by $g,x_1,x_2$.

By \cite{g}, the function $\frac{c(i,j)}{c(i-1,j)}q^j$ is
$p-$periodic in each variable.
Moreover, the coproduct of $\hat H$ maps
$x_i$ into $\hat H\otimes \hat H$;
thus, similarly to \cite{g},
the function $\frac{c(i,j)}{c(i-d_k,j)}q^{b_k j}$ is $p-$periodic
in each variable for $k=1,2$.
Hence the function $\frac{c(i,j)}{c(i-1,j)}q^{(b_k/d_k) j}$ is
$p-$periodic in each variable for $k=1,2$ (here $b_k/d_k$ is the ratio taken in
$\Z_{p^2}$). We thus conclude that $b_k=d_k$ modulo $p$,
for $k=1,2$.

Now set
$\bar{g}:=g^{b_1}$, $\bar{q}:=q^{d_1b_1}$, $b:=b_2/b_1$ and
$d:=d_2/d_1$.
We obtain
$$\bar{g}x_1\bar{g}^{-1}=\bar{q}x_1,\;\Delta(x_1)=x_1\ot \bar{g}+1\ot x_1$$
and
$$\bar{g}x_2\bar{g}^{-1}=\bar{q}^{d}x_2,\;\Delta(x_2)=x_2\ot \bar{g}^{b}
+1\ot x_2,$$
where $b,d\in \Z_{p^2}$ and $b=d$ modulo $p$.

Extend $\bar{H}$ to a Hopf algebra $H'$ generated by $\bar{H}$
and two commuting grouplike elements $g_1,g_2$, with relations
$g_ix_jg_i^{-1}=\bar{q}^{\delta_{ij}}
x_j$, $g_i^{p^2}=1$ for $i,j=1,2$, and $\bar g=g_1g_2^d$.
(The proof that this is possible is the same as
the proof given above of the fact
that $H$ can be extended by adjoining $g$.)

Let $\lambda\in \Z_{p^2}$.
Let $$T=T_\lambda:=\sum_{\gamma,\beta}\bar{q}^{\lambda\beta_1\gamma_2}
1_{\beta}\ot 1_{\gamma}\in \C[\Z_{p^2}\times
\Z_{p^2}]^{\otimes 2},$$ where $\beta=(\beta_1,\beta_2)$,
$\gamma=(\gamma_1,\gamma_2)$ and $\{1_{\beta}|\beta\in \Z_{p^2}\times \Z_{p^2}\}$ is the set
of primitive idempotents of $\Z_{p^2}\times \Z_{p^2}$.
This is a Hopf twist.
Consider the new coproduct $\Delta_T$, obtained by twisting
$\Delta$ by $T$. That is, $\Delta_T(z)=T\Delta(z)T^{-1}$.

Using the facts
that
$1_{\beta}g_i=\bar{q}^{\beta_i}1_{\beta}$ and
$1_{\beta}x_i=x_i1_{\beta-\epsilon_i}$,
$i=1,2$, where $\epsilon_1:=(1,0)$ and $\epsilon_2:=(0,1)$,
it is straightforward to verify that
$$\Delta_T(x_1)=x_1\ot g_1g_2^{\lambda+d} + 1\ot x_1\;\text{and}\;
\Delta_T(x_2)=x_2\ot g_1^bg_2^{bd} + g_1^{\lambda}\ot x_2.$$
Therefore, if we set $$z_1:=x_1,\; z_2:=g_1^{-\lambda}x_2,\;
h_1:=g_1g_2^{\lambda+d}\;
\text{and}\;h_2:=g_1^{b-\lambda}g_2^{bd}$$ we get
$$\Delta_T(z_1)=z_1\ot h_1+ 1\ot z_1\;\text{and}\;\Delta_T(z_2)=z_2\ot h_2 + 1\ot z_2.$$
Now, the relations
$$h_1z_1h_1^{-1}=\bar{q}z_1,\;h_1z_2h_1^{-1}=\bar{q}^{\lambda +d}z_2,\;
h_2z_1h_2^{-1}=\bar{q}^{b-\lambda}z_1\;\text{and}\;h_2z_2h_2^{-1}=\bar{q}^{bd}z_2,$$ imply that the
braiding matrix $B$ of $(H')^T$ (in the sense of \cite{as})
is given by $b_{11}=\bar{q}$, $b_{12}=\bar{q}^{\lambda
+d}$, $b_{21}=\bar{q}^{b-\lambda}$
and $b_{22}=\bar{q}^{bd}$.

Now set $\lambda=(b-d)/2$. In this case
$b_{12}=b_{21}=\bar{q}^{(b+d)/2}$, so the
braiding matrix is symmetric, and the corresponding Nichols
algebra is of FL type in the sense of \cite{as}.

According to \cite{as}, the
Cartan matrix $A$ corresponding to $B$ has $a_{12}=b+d$
and $a_{21}=(b+d)/bd$ (modulo $p^2$).
Since in our situation $b=d$ modulo $p$, we get that modulo $p$, $a_{12}=2b$
and $a_{21}=2/b$, and hence that $a_{12}a_{21}=4$ modulo $p$.
We claim that this implies that the Cartan matrix $A$ cannot be of finite
type.

Indeed, in the finite type case
($A_1\times A_1$, $A_2$, $B_2$ and $G_2$),
$a_{12}a_{21}=0,1,2,3$. Therefore if $p>3$, $A$ cannot be of
finite type.
For $p=3$, we get that $a_{12}=a_{21}=-1$ ($A_2$ case) and
$b=1$ modulo $3$. But this implies
that $b^2+b+1=0$ modulo $9$, which leads to a contradiction.

Now by Theorem 1.1 (ii) in \cite{as}, the algebra
$(H')^T$ (and hence $H'$) is infinite dimensional.
This gives a contradiction and completes the proof of Theorem \ref{rank}.

\section{Construction of finite dimensional basic quasi-Hopf algebras}

In this section we generalize the construction of $A(q)$ from \cite{g}, and
construct finite dimensional basic quasi-Hopf algebras which are not twist
equivalent to a Hopf algebra.

Let $n\ge 2$ be an integer, and $q$ a primitive root of $1$ of order
$n^2$. Let $H$ be a finite dimensional Hopf algebra generated by
grouplike elements $g_i$ and skew-primitive elements $e_i$, $i=1,\dots,m$,
such that
$$g_i^{n^2}=1,\;g_ig_j=g_jg_i,\;g_ie_jg_i^{-1}=q^{\delta_{i,j}}e_j$$
and
$$\Delta(e_i)=e_i\ot K_i+1\ot e_i,$$
where $K_i:=\prod_jg_j^{a_{ij}}$ for some $a_{ij}$ in $\Z_{n^2}$.

Assume that $H$ has a projection onto $\C[(\Z_{n^2}) ^m]$, $g_i\mapsto g_i$
and $e_i\mapsto 0$, and let $B\subset H$ be the subalgebra generated by $\{e_i\}$.
Then by Radford's theorem \cite{r}, the multiplication map $\C[(\Z_{n^2}) ^m]\ot
B\to H$ is an isomorphism of vector spaces. Therefore, $A:=\C[(\Z_n) ^m]B\subset H$
is a subalgebra of dimension $\dim(H)/n^m$.
It is generated by $g_i^n$ and $e_i$.

Let $\{1_{\beta}|\beta=(\beta_1,\dots,\beta_m)\in (\Z_{n^2}) ^m\}$ be the set
of primitive idempotents of $\Bbb C[\Z_{n^2}^m$, and denote by $\epsilon_i\in (\Z_{n^2}) ^n$ the vector
with $1$ in the $i$th place and $0$ elsewhere. Note that
$$1_{\beta}g_i=q^{\beta_i}1_{\beta}\;\text{and}\;1_{\beta}e_i=e_i1_{\beta -\epsilon_i}.$$

Let $c(z,y)$ be the coefficients of the twist $J$ as above
introduced in \cite{g}. Recall from \cite{g} that
$c(z,y)=q^{-z(y-y')}$, where $y'$ denotes
the remainder of division of $y$ by $n$.

Let
$${\mathbb J}:=\sum_{\beta,\gamma\in (\Z_{n^2}) ^m}\prod_{i,j=1}^{m}c(\beta_i,\gamma_j)
^{a_{ij}}1_{\beta}\ot 1_{\gamma}.$$
It is clear that it is invertible and
$(\varepsilon\ot \id)({\mathbb J})=(\id\ot \varepsilon)({\mathbb J})=1$.
Define a new coproduct $\Delta_{\mathbb J}(z)=\mathbb J\Delta(z)\mathbb
J^{-1}$.

\begin{lemma} The elements
$\Delta_{\mathbb J}(e_i)$ belong to $A\ot A$.
\end{lemma}

\begin{proof}
This lemma for $m=1$ was proved in \cite{g}.
The general case follows from the case $m=1$
by a straightforward computation.
\end{proof}

\begin{lemma}
The associator $\Phi:=d{\mathbb J}$
obtained by twisting the trivial associator by ${\mathbb J}$
is given by the formula
$$\Phi=\sum_{\beta,\gamma,\delta\in \Z_n^m}\biggl(\prod_{i,j=1}^m
q^{a_{ij}\beta_i((\gamma_j+\delta_j)'-\gamma_j-\delta_j)}\biggr)\mathbf
1_\beta\otimes \mathbf 1_\gamma\otimes \mathbf 1_\delta,$$
where $\{\mathbf 1_\beta\}$ are the primitive idempotents of
$\Z_n^m$ ($\mathbf 1_\beta g_i^n=q^{n\beta_i}\mathbf 1_\beta$),
and we regard the components of $\beta,\gamma,\delta$
as elements of $\mathbb Z$.
Thus $\Phi$ belongs to $A\otimes A\otimes A$.
\end{lemma}

\begin{proof} One has
$$
\Phi=\sum_{\beta,\gamma,\delta\in (\Z_{n^2}) ^m}
\prod_{i,j=1}^m\left (\frac{c(\beta_i+\gamma_i,\delta_j)c(\beta_i,\gamma_j)}
{c(\beta_i,\gamma_j+\delta_j)c(\gamma_i,\delta_j)}\right )
^{a_{ij}}1_{\beta}\ot 1_{\gamma}\ot 1_{\delta}.
$$
Substituting the expression of $c(z,y)$,
similarly to \cite{g} we get the statement.
\end{proof}

Thus we get our second main result.

\begin{theorem} The algebra $A$ is a quasi-Hopf subalgebra
of $H^{\mathbb J}$, which has coproduct $\Delta_{\mathbb J}$ and
associator $\Phi$.
\end{theorem}

\begin{proof}
We have shown that $\Delta_{\mathbb J}:A\to A\ot A$ and $\Phi\in
A\ot A\ot A$. It is also straightforward to show that $S_{\mathbb J}:A\to A$
and $\alpha\in A$ if $\beta$ is gauged to be $1$
(where $S_{\mathbb J}$, $\alpha$, and $\beta$ are the antipode and
the distinguished elements of $H^{\mathbb J}$).
Thus $A$ is a quasi-Hopf subalgebra of $H^{\mathbb J}$.
\end{proof}

This yields many examples of
finite dimensional basic quasi-Hopf algebras $A$. For instance,
let $\mathfrak{g}$ be a finite dimensional simple Lie algebra,
and $\mathfrak{b}$ be a Borel subalgebra of $\g$.
Assume that $n$ is odd, coprime to $3$ if $\g=G_2$.
Then we can take
$H$ to be the Frobenius-Lusztig kernel
$u_q(\mathfrak{b})$. In this case, $A$ is a
quasi-Hopf algebra of dimension $n^{{\rm dim}{\mathfrak g}}$.
Another example is obtained from $H=\gr(u_q(\mathfrak{g}))$ (with respect
to the coradical filtration).

\begin{remark} If for some $i$, $a_{ii}\ne 0$ modulo $n$,
then $A$ is not twist equivalent to a Hopf algebra.
Indeed, the associator $\Phi$ is non-trivial since the
$3-$cocycle corresponding to $\Phi$ restricts to a
non-trivial $3-$cocycle on the cyclic group $\Z_n$
consisting of all tuples whose coordinates equal $0$, except for
the $i$th coordinate. Since $A$ projects onto $(\C[\Z_n^m],\Phi)$
with non-trivial $\Phi$, $A$ is not twist equivalent to
a Hopf algebra.

For instance, this is the case in the above two examples
obtained from $u_q(\mathfrak{b})$ and $u_q({\mathfrak{g}})$.
\end{remark}

\end{document}